\documentclass[12pt]{article}

\usepackage{epsfig}
\usepackage{psfrag}
\usepackage{amsbsy}
\usepackage{latexsym}
\usepackage[mathscr]{eucal}
\usepackage{amsfonts}
\usepackage{amssymb}
\usepackage{color}

\begin{document}

\title{Orthogonal Matching Pursuit under the Restricted Isometry Property
\thanks{
  This research was supported by   by the ONR Contracts
  N00014-11-1-0712,  N00014-12-1-0561, N00014-15-1-2181; the  NSF Grants  DMS 1222715,  
      DMS 0915231, DMS 1222390;
  the Institut Universitaire de France; the ERC Adv grant BREAD; 
    the DFG SFB-Transregio 40; 
    the DFG Research Group 1779;
   the Excellence Initiative of the German Federal and State Governments,
 and   RWTH Aachen  Distinguished Professorship, Graduate School AICES.
 }}

\author{Albert Cohen, Wolfgang Dahmen, and Ronald DeVore }

\hbadness=10000
\vbadness=10000

\newtheorem{lemma}{Lemma}[section]
\newtheorem{proposition}[lemma]{Proposition}
\newtheorem{cor}[lemma]{Corollary}
\newtheorem{theorem}[lemma]{Theorem}
\newtheorem{remark}[lemma]{Remark}
\newtheorem{example}[lemma]{Example}
\newtheorem{definition}[lemma]{Definition}
\newtheorem{proper}[lemma]{Properties}

\def\ds{\displaystyle}
\def\ev#1{\vec{#1}}      
\def\Supp#1{{\rm supp\,}{#1}}
\newcommand{\e}{\varepsilon}
 \def\R{\mathbb{R}}
\def\N{\mathbb{N}}
\def\C{\mathbb{C}}
\def\Z{\mathbb{Z}}
\def\E{{\relax\ifmmode I\!\!E\else$I\!\!E$\fi}}
\def\nl{\newline}
\def\T{{\relax\ifmmode I\!\!\hspace{-1pt}T\else$I\!\!\hspace{-1pt}T$\fi}}

\def\Zd{\Z^d}
\def\Q{{\relax\ifmmode I\!\!\!\!Q\else$I\!\!\!\!Q$\fi}}

\def\Rd{\R^d}
\def\gsim{\mathrel{\raisebox{-4pt}{$\stackrel{\textstyle>}{\sim}$}}}
\def\sime{\raisebox{0ex}{$~\stackrel{\textstyle\sim}{=}~$}}
\def\lsim{\raisebox{-1ex}{$~\stackrel{\textstyle<}{\sim}~$}}
\def\div{\mbox{ div }}
\def\M{M}  \def\NN{N}                   
\def\L{{\ell}}              
\def\Le{{\ell_1}}             
\def\Lz{{\ell_2}}
\def\Let{{\tilde\ell_1}}     
\def\Lzt{{\tilde\ell_2}}
 \def\t#1{\tilde{#1}}
\def\la{\lambda}
\def\La{\Lambda}
\def\ga{\gamma}
\def\BV{{\rm BV}}
\def\Ga{\eta}
\def\al{\alpha}
\def\cZ{{\cal Z}}
\def\pr{\bf P}
\def\argmax{\mathop{\rm Argmax}}
\def\argmin{\mathop{\rm Argmin}}
 \def\cO{{\cal O}}
\def\cA{{\cal A}}
\def\cC{{\cal C}}
\def\cF{{\cal F}}
\def\bu{{\bf u}}
\def\bz{{\bf z}}
\def\bZ{{\bf Z}}
\def\bI{{\bf I}}
\def\cE{{\cal E}}
\def\cD{{\cal D}}
\def\cG{{\cal G}}
\def\cI{{\cal I}}
\def\cJ{{\cal J}}
\def\cM{{\cal M}}
\def\cN{{\cal N}}
\def\cT{{\cal T}}
\def\cU{{\cal U}}
\def\cV{{\cal V}}
\def\cW{{\cal W}}
\def\cL{{\cal L}}
\def\cB{{\cal B}}
\def\cG{{\cal G}}
\def\cK{{\cal K}}
\def\cS{{\cal S}}
\def\cP{{\cal P}}
\def\cQ{{\cal Q}}
\def\cR{{\cal R}}
\def\cU{{\cal U}}
\def\bL{{\bf L}}
\def\bK{{\bf K}}
\def\bC{{\bf C}}
\def\X{X\in\{L,R\}}
\def\ph{{\varphi}}
\def\D{{\Delta}}
\def\H{{\cal H}}
\def\bM{{\bf M}}
\def\bG{{\bf G}}
\def\bP{{\bf P}}
\def\bW{{\bf W}}
\def\bT{{\bf T}}
\def\bV{{\bf V}}
\def\bv{{\bf v}}
\def\bz{{\bf z}}
\def\bw{{\bf w}}
\def \span{{\rm span}}
\def \meas {{\rm meas}}
\def\rhom{{\rho^m}}
 \def\lll{\langle}
\def\rr{\rangle}
\def\dJ{\nabla}
\newcommand{\ba}{{\bf a}}
\newcommand{\bb}{{\bf b}}
\newcommand{\bc}{{\bf c}}
\newcommand{\bd}{{\bf d}}
\newcommand{\bs}{{\bf s}}
\newcommand{\bff}{{\bf f}}
\newcommand{\bp}{{\bf p}}
\newcommand{\bg}{{\bf g}}
\newcommand{\sm}{\setminus}
\newcommand{\br}{{\bf r}}
 
\def \Vol{\mathop{\rm  Vol}}
\def \mes{\mathop{\rm mes}}
\def \Prob{\mathop{\rm  Prob}}
\def \exp{\mathop{\rm    exp}}
\def \sign{\mathop{\rm   sign}}
\def \sp{\mathop{\rm   span}}
\def \vphi{{\varphi}}
\def \csp{\overline \mathop{\rm   span}}
\def\o{\overline}

\def\<{\langle}
\def\>{\rangle}
\def\[{\Bigl [}
\def\]{\Bigr ]}
\def\({\Bigl (}
\def\){\Bigr )}
\def\ul#1{\underline{#1}}

\newenvironment{disarray}{\everymath{\displaystyle\everymath{}}\array}{\endarray}

\newcommand{\be}{\begin{equation}}
\newcommand{\ee}{\end{equation}}
\newcommand{\beqn}{\begin{equation}}
\newcommand{\eeqn}{\end{equation}}
\def\beginproof{\noindent{\bf Proof:}~ }
\def\endproof{\hfill\rule{1.5mm}{1.5mm}\\[2mm]}

 \newenvironment{Proof}{\noindent{\bf Proof:}\quad}{\endproof}

\renewcommand{\theequation}{\thesection.\arabic{equation}}
\renewcommand{\thefigure}{\thesection.\arabic{figure}}

\makeatletter
\@addtoreset{equation}{section}
\makeatother

\newcommand\abs[1]{\left|#1\right|}
\newcommand\clos{\mathop{\rm clos}\nolimits}
\newcommand\trunc{\mathop{\rm trunc}\nolimits}
\renewcommand\d{d}
\newcommand\dd{d}
\newcommand\diag{\mathop{\rm diag}}
\newcommand\dist{\mathop{\rm dist}}
\newcommand\diam{\mathop{\rm diam}}
\newcommand\cond{\mathop{\rm cond}\nolimits}
\newcommand\eref[1]{(\ref{#1})}
\newcommand\Hnorm[1]{\norm{#1}_{H^s([0,1])}}
\def\int{\intop\limits}
\renewcommand\labelenumi{(\roman{enumi})}
\newcommand\lnorm[1]{\norm{#1}_{\ell_2(\Z)}}
\newcommand\Lnorm[1]{\norm{#1}_{L_2([0,1])}}
\newcommand\LR{{L_2(\R)}}
\newcommand\LRnorm[1]{\norm{#1}_\LR}
\newcommand\Matrix[2]{\hphantom{#1}_#2#1}
\newcommand\norm[1]{\left\|#1\right\|}
\newcommand\ogauss[1]{\left\lceil#1\right\rceil}
\newcommand{\QED}{\hfill
\raisebox{-2pt}{\rule{5.6pt}{8pt}\rule{4pt}{0pt}}
  \smallskip\par}
\newcommand\Rscalar[1]{\scalar{#1}_\R}
\newcommand\scalar[1]{\left(#1\right)}
\newcommand\Scalar[1]{\scalar{#1}_{[0,1]}}
\newcommand\Span{\mathop{\rm span}}
\newcommand\supp{\mathop{\rm supp}}
\newcommand\ugauss[1]{\left\lfloor#1\right\rfloor}
\newcommand\with{\, : \,}

\newcommand\Null{{\bf 0}}
\newcommand\bA{{\bf A}}
\newcommand\bB{{\bf B}}
\newcommand\bR{{\bf R}}
\newcommand\bD{{\bf D}}
\newcommand\bE{{\bf E}}
\newcommand\bF{{\bf F}}
\newcommand\bH{{\bf H}}
\newcommand\bU{{\bf U}}
\newcommand\ve{\varepsilon}
\newcommand\cH{{\cal H}}
\newcommand\sinc{{\rm sinc}}
\def\enorm#1{| \! | \! | #1 | \! | \! |}

\newcommand{\dm}{\frac{d-1}{d}}

\let\bm\bf
\newcommand{\bbeta}{{\mbox{\boldmath$\beta$}}}
\newcommand{\bal}{{\mbox{\boldmath$\alpha$}}}
\newcommand{\bbi}{{\bm i}}

\newif\ifNZB
\newcommand\NZB[1]{\ifNZB \marginpar{\raggedright \scriptsize NZB:\\#1}
 \else \fi}
\newcommand{\FText}[1]{\mbox{#1}}
\makeatletter
\newcommand{\tr}{{\mathop{\operator@font T}\nolimits}}
\newcommand{\mod}{\mathop{\operator@font mod}}
\makeatother
\newcommand{\UArrow}[4]{
 \begin{array}{ll}
  #1&\stackrel{#2}\longrightarrow\\
  #3&\;\raisebox{1ex}{$\nearrow$}\mkern-14mu_{#4}
 \end{array}}
\newcommand{\DArrow}[4]{
 \begin{array}{ll}
  \stackrel{#1}\longrightarrow&#2\\
  \mkern-10mu_{#3}\mkern-22mu\raisebox{1ex}{$\searrow$}&#4
 \end{array}}
\newcommand{\fig}[3]{\par\begin{figure}[ht]
  \centerline{\epsfbox{#1.eps}}\caption{#3}\label{fig#2}\end{figure}}
 \newcommand{\dI}{\Delta}
\newcommand{\Cn}{\bar C_0}

\newcommand{\Phix}{\Phi x}
\newcommand{\PhixT}{\Phi x_T}
\newcommand{\Phiet}{\Phi\eta}
\newcommand{\Pho}{\Phi x_0}
\newcommand{\Phu}{\Phi x_1}
\newcommand{\Phd}{\Phi x_2}
\newcommand{\Phiz}{\Phi z}
\newcommand{\Phietu}{\Phi\eta_1}
\newcommand{\Phieto}{\Phi\eta_0}
\newcommand{\vp}{\varphi}

\maketitle
\date{}

\begin{abstract}
  This paper is concerned with the performance of Orthogonal Matching Pursuit (OMP) algorithms 
  applied to a dictionary $\cD$ in a Hilbert space $\cH$.    Given an element $f\in \cH$, OMP generates 
  a sequence of approximations $f_n$, $n=1,2,\dots$, each of which is a linear combination of $n$ 
  dictionary elements chosen by a greedy criterion.  It is studied whether the approximations $f_n$ 
  are in  some  sense comparable to {\em best $n$ term approximation} from the dictionary.
  One important result related to this question is a theorem of   
  Zhang \cite{TZ} in the context of sparse recovery of finite dimensional signals. This
  theorem shows that OMP exactly recovers 
  $n$-sparse signal, whenever the dictionary $\cD$ satisfies a Restricted Isometry Property (RIP)
  of order $An$ for some constant $A$, and that the procedure is
  also stable in $\ell^2$ under measurement noise. 
The main contribution of the present paper is to give a structurally simpler proof of
  Zhang's theorem, formulated in the general context of $n$ term approximation 
  from a dictionary in arbitrary Hilbert spaces $\cH$. Namely, it is shown that OMP
  generates near best $n$ term approximations under a similar RIP condition.  
  
  \end{abstract}


\section{Introduction}
\label{intro}

Approximation by  sparse linear combinations of elements from a fixed redundant family is a frequently
employed technique in signal processing and other application domains.   
 We consider such problems in
  a separable Hilbert space $\cH$ endowed with a norm $\|\cdot\|:=\|\cdot\|_\cH$ induced by the scalar product $\lll\cdot,\cdot\rr$
on $\cH\times \cH$. A countable collection $\cD= \{\vp_\gamma\}_{\gamma\in\Gamma}\subset \cH$ is called a {\it dictionary} if it is complete, i.e.,
the set of finite linear combinations of elements of the dictionary are dense in $\cH$.   The simplest example of a dictionary is the set of elements of a fixed basis of $\cH$.   But our primary interest is in {\em redundant} families.   In such a case, there exists a strict subset of $\cD$ that is still a dictionary.  A primary example  of a redundant dictionary  is a frame, e.g.,  any union of a finite number of bases.   
Without loss of generality we shall always
assume that the dictionary $\cD$ is {\em normalized}, i.e., 
$$
\|\vp_\gamma\| =1,\quad \gamma \in \Gamma.
$$

Given such a dictionary $\cD$, we consider the class
\be
\Sigma_n=\Sigma_n(\cD):=\Big\{\sum_{\gamma\in S} c_\gamma\vp_\gamma\; : \; \#(S)\leq n\Big\} \subset \cH,\quad n\ge 1.
\label{13Sigman}
\ee
 The elements in $\Sigma_n$  
are said to be  sparse with {\em sparsity} $n$.            We define
$$
\sigma_n(f)_\cH := \inf_{g\in \Sigma_n}\|f-g\|,
$$ 
which is called {\it the error of best $n$-term approximation} to
$f$ from the dictionary $\cD$.    

An important distinction between  $n$ term  dictionary approximation and   other forms of approximation, such as approximation from an $n$ dimensional space, is that the set $\Sigma_n$ is not a linear space since the sum of two elements in $\Sigma_n$ is generally not in $\Sigma_n$, although it is in $\Sigma_{2n}$. 
Thus   $n$-term approximation from a dictionary is  an important example of  nonlinear approximation  \cite{DActa}  that 
reaches into numerous application areas such as adaptive PDE solvers, image encoding, or statistical learning. 
It also serves as a performance benchmark in {\em compressed
sensing} that better captures the robustness of compressed sensing than results on exact sparsity recovery \cite{CDD-cs}.

While there are many themes in $n$ term dictionary approximation, our interest here is in analyzing the performance of   greedy algorithms  for generating $n$-term approximations to a given target element $f\in \cH$.    There are numerous papers on this subject.  We refer the reader to the survey
article \cite{TemActa} as a general reference.   Our particular interest is in understanding what properties of the dictionary $\cD$ guarantee that these algorithms
perform similarly to best $n$-term approximation.

 These algorithms and best $n$-term approximation have a simple description when the dictionary $\cD$ is  an orthonormal or, more generally, a Riesz basis of  $\cH$.  In this case, the  best $n$-term approximations to a given $f\in \cH$  are    realized by expanding $f$ in terms of the basis 
 \be
f=\sum_{\gamma\in\Gamma}c_\gamma\vp_\gamma
\label{13expan}
\ee
and  retaining  $n$ terms from this expansion which correspond to the  largest (in absolute value) expansion coefficients.  The typcial greedy algorithm will construct the same approximations.
The situation is much less clear when dealing with more general dictionaries.

  In the case of general dictionaries,  algorithms for generating $n$-term approximations are typically built
  on some form of greedy selection
  \be
\vp_k:=\vp_{\gamma_k}, \;\; k=1,2,\dots,
\ee
of elements from $\cD$ and then using a linear combination of $ \vp_1,\dots,\vp_n$ as the $n$-term approximation.  
   The standard greedy algorithm (called the Pure Greedy Algorithm) makes the initial selection $\vp_1$  as any element  such that
\be
\label{initialomp}
\vp_1:=\argmax_{\vp\in\cD} | \langle f,\vp\rangle|.
\ee
This gives the approximation $f_1:=\langle f,\vp_1\rangle \vp_1$ to $f$ and the residual $r_1:=f-f_1$.
Given that $\vp_1,\dots,\vp_{k-1}$ have been selected, and an approximation $f_{k-1}$ from $F_{k-1}:=\span\{\vp_1,\dots,\vp_{k-1}\}$ has been constructed, 
 the next dictionary element 
$\vp_k$ is chosen as the  best match of the residual 
\be
\label{res1}
r_{k-1}:=f-f_{k-1},
\ee
in the sense that
\be
\vp_{k}:=\argmax_{\gamma\in\Gamma} |\<r_{k-1},\vp_\gamma\>|.
\label{13matching}
\ee

There exist different ways of forming  the next approximation $f_k$ resulting in different greedy  algorithms.   We focus our attention on    {\em Orthogonal Matching Pursuit }(OMP),
  which forms  the new approximation as 
\be
f_k:=P_k f,
\ee
where $P_k$ is the orthogonal projector onto $F_k$.  OMP is also called the Orthogonal Greedy Algorithm.
More generally, we analyze the   {\em Weak Orthogonal Matching Pursuit} (WOMP) where the choice of 
$\vp_k$ is only required to satisfy 
\be
|\<r_{k-1},\vp_k\>|\geq \kappa \max_{\gamma\in\Gamma} |\<r_{k-1},\vp_\gamma\>|,
\ee
where $\kappa\in ]0,1]$ is a fixed parameter, which is a more easily implemented selection rule in
practical applications.   Once this choice of $\vp_1,\dots,\vp_k$ is made, then $f_k$ is again defined as the orthogonal projection onto $F_k$.

The main interest of the present paper is to understand what properties of a dictionary $\cD$ guarantee that  the approximation rate of WOMP
after $O(n)$ steps is comparable to the   the best $n$-term approximation error  $\sigma_n(f)$, at least for a certain range $n\leq N$. A  related question,  but less demanding, is to understand when
  WOMP is  guaranteed to exactly recover $f$ whenever $f\in \Sigma_n$ in $O(n)$ steps for a suitable range of $n$.  This is 
  sometimes refered to as {\it sparse recovery}. 
  Of course, as already mentioned, we know that both of these questions have a positive answer for the entire range of $n$ whenever $\cD$ is a Riesz basis for $\cH$.

To give a precise formulation of the type of performance we seek,  we define the concept of {\it instance optimality}.
\nl

\noindent
{\bf Instance Optimality:}  {\it We say that the WOMP algorithm satisfies   instance optimality for $n\le N$, if there are constants  $A,C>0$, with $A$ an integer,  such that the outputs $f_n$  of WOMP satisfy 
\be
\label{goal}
\|f-f_{An}\|\leq C\sigma_n(f)_\cH,
\ee
for   $n\le N$.}
\nl

 Notice that if \eref{goal} is satisfied then it implies a positive solution to the sparse recovery problem for the same range of $n$ since $\sigma_n(f)=0$ when $f$ is in $\Sigma_n$.
  To obtain results on sparse recovery or instance optimality requires structure on the dictionary $\cD$.   The first results of this type were obtained under assumptions on the
   {\em coherence} of a dictionary $\cD\subset \cH$ defined by  
$$
\mu=\mu(\cD):=\sup\{|\<\vp,\psi\>|\;  : \; \vp,\psi\in\cD,\; \vp \neq \psi \}.
$$ 
The first results on this general circle of problems centered on sparse recovery.   Tropp \cite{Tro} proved that whenever the dictionary has coherence $\mu< \frac 1 {2n-1}$,  then $n$ steps of OMP recover any $f\in \Sigma_n$ exactly.   

Concerning instance optimality, we mention that  Livschitz \cite{Liv} proved that whenever $\mu\leq \frac 1 {20n}$, 
then after $2n$ steps, the OMP algorithm returns $f_{2n}\in\Sigma_{2n}$ such that
\be
\|f-f_{2n}\|\leq 3\sigma_n(f)_\cH.
\ee
 
A weaker assumption on a dictionary, known as the {\it Restricted Isometry Property} (RIP), was introduced in 
the context of compressed sensing \cite{CRT}. To formulate this property, we introduce the notation
\be
\Phi\bc=\sum_{\gamma\in\Gamma} c_\gamma \vp_\gamma,
\ee
whenever $\bc=(c_\gamma)_{\gamma\in \Gamma}$ is a finitely supported sequence.
The dictionary $\cD$ is said to satisfy  the RIP  of order $n\in \N$ with constant   $0<\delta<1$ provided
\be
(1-\delta)\|\bc\|_{\ell^2}^2\leq \|\Phi\bc\|^2 \leq (1+\delta)\|\bc\|_{\ell^2}^2, \quad
Ê\|\bc\|_{\ell^0}:=  \#({\rm supp}\,\bc)\leq n.
\label{13rip}
\ee
Hence this property quantifies the deviation of any subset of cardinality at most $n$ from
an orthonormal set.
We denote by $\delta_n$ the minimal value of $\delta$ for which this property holds and remark that trivially $\delta_n\leq \delta_{n+1}$.

It is well-known that  a   coherence 
bound
\be
\label{smallcoherence}
\mu(\cD)< (n-1)^{-1}
\ee 
implies the validity of ${\rm RIP}(n)$ for 
$\delta_n\le (n-1)\mu$, but not vice versa \cite{Tro}. 

In \cite{TZ}, Tong Zhang proved that OMP exactly recovers 
finite dimensional $n$-sparse signals, whenever the dictionary $\cD$ satisfies a Restricted Isometry Property (RIP)
of order $An$ for some constant $A$, and that the procedure is
also stable in $\ell^2$ under measurement noise. 
The main result of the present paper is the following related 
theorem on instance optimality for WOMP.

\begin{theorem}
\label{maintheorem}
Given the weakness parameter $\kappa \le 1$, there exist fixed constants $A,C,\delta^*$, such that the following holds for all $n\geq 0$:
if $\cD$ is a dictionary in a Hilbert space $\cH$ for which 
${\rm RIP((A+1)n)}$ holds with $\delta_{(A+1)n}\leq \delta^*$, then, for any
target function $f\in \cH$, the WOMP algorithm returns after $An$ steps
an approximation $f_{An}$ to $f$ that satisfies
\be
\|f-f_{An}\| \leq C\sigma_n(f)_\cH.
\label{13ompopt}
\ee
\end{theorem}

The values of $A$, $C$, $\kappa$, and $\delta^*$ for which the above result holds
are coupled. For example, it is possible to have a smaller value of $A$
at the price of a larger value of $C$ or of a smaller value of $\delta^*$. Similarly, a smaller weakness parameter $\kappa$ can be compensated
by increasing $A$.

While the theorem of \cite{TZ} is not stated in the above form, it can be used
to derive Theorem \ref{maintheorem} by interpreting the error of best $n$-term
approximation as a measurement noise. In this way, one version of the above 
result can be derived from \cite{TZ} for OMP ($\kappa=1$) with $\delta^*=\frac 1 3$ and $A=30$. 
Let us mention that Zhang's theorem is also established in \cite{FR}, with the same proof,
but with different constants $\delta^*=\frac 1 6$ and $A=12$.    

In what follows, we do not focus on improving the constants, but
rather our interest is to  provide a conceptually more elementary proof
for Theorem \ref{maintheorem}. Namely the proof for  \cite{TZ}Ê and \cite{FR} 
is based on an induction argument which involves an auxiliary greedy algorithm
(initialized from a non trivial sparse approximation) in an inner loop. Our proof
avoids using this auxiliary step. It is also presented in the framework of a
possibly infinite dimensional Hilbert space $\cH$. We give the new proof in the following section. 
We then give some observations that can be derived from  
Theorem \ref{maintheorem}.

In this paper, we shall sometimes use the notation 
$\Phi^* v=(\<v,\vp_\gamma\>)_{\gamma\in \Gamma}$ for any $v\in \cH$,
and $\bc_T$ to denote, for any $\bc=(c_\gamma)_{\gamma\in \Gamma}$ and $T\subset \Gamma$,
the sequence whose entries coincides with those of 
$\bc$ on $T$ and are $0$ otherwise.

\section{Proof of Theorem \ref{maintheorem}}\label{sec:3}

In this section, we give a proof of Theorem \ref{maintheorem}.  We begin with the following elementary lemma
which guarantees the existence of near best $n$ term approximations from a dictionary.

\begin{lemma}
\label{13lemclose}
Let $\cD$ be a dictionary in a Hilbert space $\cH$ that 
satisfies $RIP(2n)$. 
Then,

\noindent
{\rm (i)} the set $\Sigma_n$
of all $n$-term linear combinations from $\cD$ is closed in $\cH$.

\noindent
{\rm (ii)} For each $f\in \cH$, $\e>0$, and $n\ge 1$, there exists a $g\in \Sigma_n$ such that
\be
\label{nb}
\|f-g\|\le (1+\e)\sigma_n(f)_\cH.
\ee
\end{lemma}

\noindent
{\bf Proof:} To prove (i), we let $(g^k)_{k\geq 0}$ be a sequence of 
elements from $\Sigma_n$ that converges in $\cH$
towards some $g\in \cH$. We may write
\be
g^k=\Phi \bc^k=\sum_{\gamma\in\Gamma} c^k_\gamma\vp_\gamma,
\ee
with $\|\bc^k\|_{\ell^0} \leq n$. For any $\e>0$, there exists $K$ such that
\be
\|g^k-g^l\| \leq \e, \quad k,l\geq K.
\ee
From RIP$(2n)$, it follows that
\be
\|\bc^k-\bc^l\|_{\ell^2} \leq \frac{\e}  {\sqrt{1-\delta_{2n}}} ,
\ee
which shows that the sequence $(\bc^k)_{k\geq 0}$ converges in $\ell^2$
to some $\bc\in \ell^2$. In particular, we find that
\be
\lim_{k\to +\infty} c^k_\gamma=c_\gamma,\quad \gamma\in \Gamma.
\ee
If $c_\gamma\neq 0$ for more than $n$ values of $\gamma$, we find
that $\|\bc^k\|_{\ell^0}>n$  for $k$ sufficiently large which is a contradiction.
It follows that $g=\sum_{\gamma\in\Gamma}c_\gamma\vp_\gamma\in\Sigma_n$.

To prove (ii), let $g^k\in\Sigma_n$  be such that $\|g_k-f\|\to \sigma_n(f)_\cH$.  If $\sigma_n(f)>0$, then $g=g_k$ will satisfy (ii) if $k$ is sufficiently large.   On the other hand, if $\sigma_n(f)=0$, then $g_k\to f$, $k\to\infty$.    By (i) $f\in \Sigma_n$ and so we can take $g=f$.\hfill $\Box$

\subsection{Reduction of the residual}

Our  starting point in proving Theorem \ref{maintheorem} is the following lemma from \cite{TZ} which quantifies  the reduction of the residuals generated by the WOMP algorithm under the RIP condition. 
In what follows, we denote by
\be
S_k:=\{\gamma_1,\dots,\gamma_k\},
\ee
the set of indices selected after $k$ steps of WOMP applied to
the given target element $f\in \cH$, and denote as before the residual by $r_k=f-f_k$. 
\begin{lemma}
\label{TZLemma}
Let $(f_k)_{k\geq 0}$ be the sequence of approximations generated by the WOMP algorithm
applied to   $f$, and let $g=\Phi \bz$ with $\bz$ supported on 
a finite set $T$. Then, if $T$ is not contained in $S_k$, one has
\be
\label{13lemdec1}
\|r_{k+1}\|^2\leq \|r_k\|^2- \frac {\kappa^2(1-\delta)}{\#(T\sm S_k)}\max\{0, \|r_k\|^2-\|f-g\|^2\},
\ee
where $\delta:=\delta_{\#(T\cup S_k)}$ is the corresponding RIP-constant and $\kappa\in ]0,1[$ is the weakness parameter in the
WOMP algorithm.
\label{13lemdec}
\end{lemma}

For completeness, we recall the proof at the end of this section.
 It is at this point, we depart from the arguments in \cite{TZ} with the goal of providing a simpler more
 transparent argument.    
 An immediate consequence of Lemma \ref{TZLemma} is the following. 
\begin{proposition}
\label{13prop:red}
Assume that for a given $A >0$ and  $\delta^*<1$, 
${\rm RIP((A+1)n)}$ holds with $\delta_{(A+1)n}\leq\delta^*$.
If $g=\Phi\bz$, where $\bz$ is supported on a set $T$
such that $\#(T)\leq n$, then for any non-negative integers $(j,m,L)$
such that $\#(T\setminus S_j)\leq m$ and 
$j+mL \leq An$, one has
\be
\|r_{j+mL}\|^2\leq e^{-\kappa^2(1-\delta^*)L}\|r_{j}\|^2+\|f-g\|^2.
\label{13iterate}
\ee
\end{proposition}

\noindent
{\bf Proof:} By Lemma \ref{13lemdec},
if $g=\Phi\bz$ where $\bz$ is supported on a set $T$
such that $\#(T)\leq n$, then for any non-negative integers $(j,m,L)$
such that $\#(T\setminus S_j)\leq m$ and 
$j+mL \leq An$, one has
$$
\begin{disarray}{ll}
\max\{0,\|r_{j+mL}\|^2-\|f-g\|^2\}& \leq \(1- \kappa^2(1-\delta^*)/ m\)^{mL} 
\max\{0,\|r_{j}\|^2-\|f-g\|^2\}\\
&\leq e^{-\kappa^2(1-\delta^*)L}\max\{0,\|r_{j}\|^2-\|f-g\|^2\}\,,
\end{disarray}
$$
where we have used the fact that $\#(T\setminus S_l)\le m$ for all $l\ge j$,
This gives \eref{13iterate} and completes the proof of Proposition \ref{13prop:red}. \hfill $\Box$
 \nl
 \nl
{\bf Proof of Theorem \ref{maintheorem}:} We fix $f$ and  use  
the abbreviated notation
\be
\label{13abbreve}
\sigma_n:=\sigma_n(f)_\cH,\quad n\ge 0.
\ee
We first observe that the assertion of the theorem follows from the following.
\nl

\noindent
{\bf Claim:} {\it  If $0\le k<n$ satisfies
\be
\label{13rAk}
\|r_{Ak}\|\le 2\sigma_k,
\ee
and is such that $\sigma_n< \frac {\sigma_k}4$, then there exists $k<k'\le n$ such that}
\be
\|r_{Ak'}\|\le 2\sigma_{k'}.
\label{13Akprime}
\ee

Indeed, assuming that this claim holds, we complete the proof of the Theorem as follows. We let $k$ be the largest integer in $\{0,\dots,n\}$
for which $\|r_{Ak}\|\le 2\sigma_k$.  Since $\|r_{0}\|=\sigma_0=\|f\|$, 
such a $k$ exists.  If $k<n$, then we must have $\sigma_k\le 4 \sigma_n$
and therefore 
\be
\|r_{An}\|\le \|r_{Ak}\|\le  2 \sigma_k\le 8\sigma_n,
\ee
so that \eref{13ompopt}
holds with $C=8$.

We are therefore left with proving the claim.  For this, we fix
\be
\label{parsetting}
\delta^*=\frac 1 6,
\ee
and $0\le k<n$   such 
that \eref{13rAk} holds and such that $\sigma_n< \frac {\sigma_k}4$.  
Let $k<K\leq n$ be the first integer such that $\sigma_{K}<\frac {\sigma_k}4$.
By (ii) of Lemma \ref{13lemclose}  
we know that for any $B>1$ there is a $g\in\Sigma_K$ with $\|f-g\|\le B\sigma_K(f)$. 
Therefore, $g$ has the form
\be
g=\Phi\bz=\sum_{\gamma\in T} z_\gamma\vp_\gamma,\quad \#(T)=K.
\ee 
The significance of $K$ is that on the one hand
\be
\|f-g\|\leq B\sigma_K<\frac B 4 \sigma_k,
\ee
while on the other hand
\be
\label{13K-1}
\sigma_k\le 4 \sigma_{K-1}.
\ee
 
To eventually apply Proposition \ref{13prop:red} for the above $g$ and  $j=Ak$, we need to bound $\#(T\sm S_{Ak})$
with $A$ yet to be specified. To this end, 
we write $K=k+M$, with $M> 0$, and observe that
if $S\subset T$ is any set with $\#(S)=M$ and $g_S:=\sum_{\gamma\in S} z_\gamma\vp_\gamma$, then
\be
\|g_S\|\ge  \|f-(g-g_S)\|- \|f-g\|
\geq  \sigma_k-B\sigma_K\ge \(1-\frac B 4\) \sigma_k,
\ee
where we have used the fact that $g-g_S\in \Sigma_k$.
Using RIP, we obtain the following lower bound for 
the coefficients of $g$:  for any set $S\subset T$ of cardinality $M$
\be
\label{13boundcoeff}
\(1-\frac B 4\)^2\sigma_k^2  \le \|g_S\|^2\le (1+\delta^*)\sum_{\gamma\in S} |z_\gamma|^2= \frac{7}{6}\sum_{\gamma\in S} |z_\gamma|^2.
  \ee
Taking for $S$ the set $S_g$ of the $M$ smallest coefficients of $g$ and noting that then for  
any  more general $S\subset T$ with $\#(S)\ge  M$, one has
$\Big(\sum_{\gamma\in S} |z_\gamma|^2\Big)/\Big(\sum_{\gamma\in S_g} |z_\gamma|^2\Big)\ge \#(S)/M$, and hence  
\be
\label{13boundcoeff2}
\frac 6 7\(1-\frac B 4\)^2\frac{ \#(S)}{M}  \sigma_k^2   \le  \sum_{\gamma\in S} |z_\gamma|^2.
\ee
For the particular set $S:=T\sm S_{Ak}$, if $\#(S)\ge M$, the above bound combined with the RIP implies
$$
\begin{array}{ll}
(1-\delta^*)\frac 6 7\(1-\frac B 4\)^2\frac{ \#(S)}{M}  \sigma_k^2 &  \le \|g_S\|^2  \leq  \|g-f_{Ak}\|^2 
 \leq (\|g-f\|+\|r_{Ak}\|)^2\\
&\leq (B\sigma_K+2\sigma_k)^2
\leq \(\frac B 4+2\)^2\sigma_k^2.
\end{array}
$$
Since $\delta^*=1/6$ this gives the bound
\be
\label{13M}
\#(T\sm S_{Ak}) \leq \frac 7 5\frac {\(\frac B 4+2\)^2} {\(1-\frac B 4\)^2} M \leq 13M,
\ee
where the second inequality is obtained by taking $B$ sufficiently close to $1$.

We proceed now verifying the claim with $k'= K-1$ when $K-1>k$ and with $k'=k+1$ otherwise.
In the first case we can use the reduction estimate provided by Proposition \ref{13prop:red}
with $j=Ak$
in combination with \eref{13K-1} to deal with the term $\|r_{Ak}\|$ in \eref{13iterate}. When $K=k+1$, however,
we cannot bound  $\|r_{Ak}\|$ directly in terms of a $\sigma_l$ for some $l>k$. Accordingly,
we use Proposition \ref{13prop:red} in different ways for the two cases.

In the case where $M\geq 2$, i.e., $K-1>k$, we apply \eref{13iterate} with $j=Ak$, $m=13M$ and $L=\lceil 4\kappa^{-2}\rceil$.
Indeed $Ak+Lm=Ak+52M\le An$ holds for $k+M\le n$ whenever $A\ge 52 \kappa^{-2}$. Moreover, notice that for  such an $A$ 
\be
A(K-1)=Ak+A(M-1)\geq Ak+\frac 1 2 AM=Ak+\frac{Am}{26} =Ak+Lm,
\ee
whenever 
\be
\label{Achoice}
A \ge 26 \lceil 4\kappa^{-2}\rceil.
\ee
This gives
$$
\begin{disarray}{ll}
\|r_{A(K-1)}\|^2 & \leq \|r_{Ak+Lm}\|^2 \\[2mm]
& \leq e^{-10/3} \|r_{Ak}\|^2+\|f-g\|^2 \\[2mm]
& \leq e^{-10/3}4\sigma_k^2+B^2\sigma_K^2\\[2mm]
& \leq e^{-10/3}64 \sigma_{K-1}^2+B^2\sigma_{K-1}^2\\[2mm]
& \leq 4\sigma_{K-1}^2,
\end{disarray}
$$
where we have used \eref{13K-1} in the fourth inequality,
and the last inequality follows by taking $B$ sufficiently close to $1$.
 We thus obtain 
\eref{13Akprime} for the value $k'=K-1>k$.

In the case $M=1$, i.e., $K=k+1$, we apply \eref{13iterate} with $j=Ak$, $m=13$ and 
$$
L=\lceil 6\kappa^{-2}\rceil.
$$
In fact, from \eref{13M} we know that $\#(T\sm S_{Ak})\le 13$ and $An \ge A(k+1)\ge Ak+mL$ for   $A$ satisfying \eref{Achoice}.
This yields
 $$
 \begin{array}{ll}
 \|r_{A(k+1)}\|^2 & \le \|r_{Ak+mL}\|^2\\[2mm]
 & \le e^{-5}\|r_{Ak}\|^2+  \|f-g\|^2 \\[2mm]
 &\le 4e^{-5}\sigma_k^2+B^2\sigma_{k+1}^2\\[2mm]
 & \le  \(4e^{-5}+\frac {B^2}{16}\)\sigma_k^2.
 \end{array}
 $$
  This implies that $S_{A(k+1)}$ contains 
 $T$. Indeed, if it missed one of the indices $\gamma\in T$, then we infer  from the RIP,
$$
\begin{array}{ll}
(1-\delta^*)|z_\gamma|^2&\le \|g-f_{A(k+1)}\|^2\\[2mm]
& \le (\|f-g\|+\|r_{A(k+1)}\|)^2\\[2mm]
&\le \(B\sigma_{K}+\sqrt{4e^{-5}+\frac {B^2}{16}}\sigma_k\)^2\\[2mm]
&\leq \(\frac B 4+\sqrt{4e^{-5}+\frac {B^2}{16}}\)^2 \sigma_k^2.
\end{array}
$$
On the other hand, we know from \eref{13boundcoeff2} that
\be
\frac 6 7\(1-\frac B 4\)^2\sigma_k^2  \le  |z_\gamma|^2,
\ee
which for $B$ sufficiently close to $1$ is a contradiction since $\frac 6 7\(1-\frac B 4\)^2> \frac 6 5\(\frac B 4+\sqrt{4e^{-5}+\frac {B^2}{16}}\)^2$.
This implies that  $\|r_{A(k+1)}\|\leq \sigma_{k+1}$, and therefore \eref{13Akprime} holds for the value $k'=k+1$.
This verifies the claim and hence completes the proof of Theorem \ref{maintheorem}.
\hfill $\Box$\\

Let us observe that Theorem \ref{maintheorem} does not give that $f_n$ is 
  a near-best $n$-term approximation in the form
\be
\label{nearbest}
\|f-f_n\| \leq C_0 \sigma_n(f)_\cH.
\ee
However a simple postprocessing of $f_{An}$ by retaining its $n$ largest components
  does satisfy \eref{nearbest}. 
\begin{cor}
\label{thm:2a}
Under the assumptions of Theorem \ref{maintheorem}, let $f_{An}=\Phi \bc^{An}$ be the output of WOMP after $An$ steps.  Let  $T\subset \Gamma$, $\#(T)=n$, be a set of indices  corresponding to $n$ largest entries of $\bc^{An}$.  Define $f_n^*\in\Sigma_n$ to be the element obtained by retaining from $f_{An}$ only the $n$ terms corresponding to the indices in $T$.  Then,
\be
\label{near-best3}
\|f-f_n^*\|\leq C^* \sigma_n(f)_\cH,
\ee
where the constant $C^*$ depends on the constant $C$ in Theorem \ref{maintheorem} and on the RIP-constant $\delta_{(A+1)n}$.
\end{cor}
{\bf Proof:}
By Lemma \ref{13lemclose}, 
  there exists a   $\bc$ with $\|\bc\|_{\ell^0}\le n$, such that
\be
\|f-\Phi \bc\| \leq 2\sigma_n(f)_\cH.
\ee
It follows that
\be
\label{13An}
\|\bc-\bc^{An}\|_{\ell^2} \leq \frac {1}{\sqrt{1-\delta_{(A+1)n}}} \|\Phi\bc-\Phi\bc^{An}\|_{\ell^2} \leq \frac {C+2}{\sqrt{1-\delta_{(A+1)n}}}\sigma_n(f)_\cH.
\ee
If   $S= {\rm supp} (\bc)$, we obtain 
\begin{eqnarray}
\|\bc-\bc^{An}_T\|_{\ell^2} & \leq & \|\bc_T-\bc^{An}_T\|_{\ell^2} + \|\bc_{T^c}-\bc^{An}_{T^c}\|_{\ell^2} + \|\bc^{An}_{T^c}\|_{\ell^2}
\nonumber\\
&\le & 2 \|\bc-\bc^{An}\|_{\ell^2} +  \|\bc^{An}_{S^c}\|_{\ell^2}\nonumber\\
&\le & 3\|\bc-\bc^{An}\|_{\ell^2},
\end{eqnarray}
which, by \eref{13An}, provides
\be
\|\bc-\bc^{An}_T\|_{\ell^2}\leq  {3}\|\bc-\bc^{An}\|_{\ell^2} \leq  \frac {  {3}(C+2)}{\sqrt{1-\delta_{(A+1)n}}}\sigma_n(f)_\cH.
\ee
The approximation $\Phi\bc^{An}_T$ is in $ \Sigma_n$ and satisfies 
\be
\|f-\Phi\bc^{An}_T\| \leq 2\sigma_n(f)_\cH+\|\Phi(\bc^{An}_T-\bc)\|
\leq \(2+\frac { {3}\sqrt {1+\delta_{(A+1)n}}(C+2)}{\sqrt{1-\delta_{(A+1)n}}}\)\sigma_n(f)_\cH,
\ee
which  proves \eref{near-best3}.\hfill $\Box$
\nl

\noindent
{\bf Proof of Lemma \ref{13lemdec}:}  
  We may assume that $\|r_k\|\geq \|f-g\|$ otherwise there is
nothing to prove. 
First observe now that
$$
\begin{disarray}{ll}
\|r_{k+1}\|^2  & =\|f-P_{k+1}f\|^2  \\[2mm]
& =\|f-P_{k}f\|^2  -\|(P_{k}-P_{k+1})f\|^2  \\[2mm]
& \leq \|r_{k}\|^2 - |\<r_{k},\vp_{\gamma_{k+1}}\>|^2.
\end{disarray}
$$
Therefore, it suffices to prove that $\|r_k\|^2-|\<r_k,\vp_{\gamma_{k+1}}\>|^2$ is bounded by the right hand side of 
\eref{13lemdec1} which amounts to showing that
\be
 (1-\delta)(\|r_k\|^2-\|f-g\|^2)\leq \kappa^{-2}\#(T\sm S_k)|\<r_k,\vp_{\gamma_{k+1}}\>|^2.
 \label{13claimTS}
 \ee
To prove this,  we first note that
$$
\begin{disarray}{ll}
2\|g-f_k\| \sqrt {\|r_k\|^2-\|f-g\|^2}& \leq  \|g-f_k\| ^2+\|r_k\|^2-\|f-g\|^2\\[2mm]
& =\ \|g-f_k\| ^2+\|r_k\|^2-\|g-f_k-r_k\|^2 \\[2mm]
&\leq 2 |\<g-f_k,r_k\>|=2|\<g,r_k\>|.
\end{disarray}
$$
This is the same as
\be
\|r_k\|^2-\|f-g\|^2\leq \frac{|\<g,r_k\>|^2}{\|g-f_k\|^2}.
\ee
If we write
$f_k=\Phi \bc^k$, with $\bc^k$   supported on $S_k$, then  the numerator of the right side satisfies
$$
\begin{disarray}{ll}
|\<g,r_k\>|&= |\<\Phi \bz,r_k\>| \\[2mm]
&=|\<\bz_{S_k^c},\Phi^*r_k\>_{\ell^2}|\\[2mm]
& \leq \|\bz_{S_k^c}\|_{\ell^1} \|\Phi^*r_k\|_{\ell^\infty}\\[2mm]
&\le \kappa^{-1}\|\bz_{S_k^c}\|_{\ell^1}|\<r_k,\vp_{\gamma_{k+1}}\>|\\[2mm]
&\leq  \kappa^{-1}\sqrt { \#(T\sm S_k)}\|\bz_{S_k^c}\|_{\ell^2}|\<r_k,\vp_{\gamma_{k+1}}\>|\\[2mm]
&  {\le}\kappa^{-1} \sqrt { \#(T\sm S_k)}\|\bz-\bc^k\|_{\ell^2}|\<r_k,\vp_{\gamma_{k+1}}\>|.
\end{disarray}
$$
On the other hand, recalling that $\delta=\delta_{\#(S_k\cup T)}$, the denominator satisfies by the RIP,
\be
\|g-f_k\|^2=\|\Phi(\bz-\bc^k)\|^2 \geq (1-\delta)\|\bz-\bc^k\|_{\ell^2}^2.
\ee
Therefore we have obtained
\be
\|r_k\|^2-\|f-g\|^2\leq \frac { \#(T\sm S_k)|\<r_k,\vp_{\gamma_{k+1}}\>|^2}{\kappa^2(1-\delta)},
\ee
which is \eref{13claimTS}.Ê\hfill $\Box$\\

 \noindent
 Albert Cohen\\ 
 Laboratoire Jacques-Louis Lions, Univerisit\'e Pierre et Marie Curie, 75005, Paris, France \\
cohen@ann.jussieu.fr
\vskip .1in
\noindent
  Wolfgang Dahmen\\ 
   Institut f\"ur Geometrie  und Praktische
Mathematik, RWTH Aachen, Templergraben 55,
D-52056 Aachen Germany\\ 
dahmen@igpm.rwth-aachen.de

 \vskip .1in
\noindent
Ronald DeVore\\
Department of Mathematics, Texas A\&M University,
College Station, T\cH 77840, USA\\
  rdevore@math.tamu.edu

\end{document}